# A NEW TYPE OF WEAKLY COMMUTATIVE GROUPS

Hariwan Z. Ibrahim [a,*], and Muwafaq M. Salih [b]

[a] Dept. of Mathematics, Faculty of Science, University of Zakho, Zakho, Kurdistan Region, Iraq - hariwan.ibrahim@uoz.edu.krd
[b] Dept. of Mathematics, College of Basic Education, University of Duhok, Kurdistan Region, Iraq - muwafaq.mahdi@uod.ac

**ABSTRACT:**
The aim of the present paper is to define and study a new class of groups, namely Wm-groups with a single binary operation based on axioms of semi commutativity, right identity and left inverse. Moreover, we introduce the notions of right cosets, quotient Wm-groups, homomorphisms, kernel and normal Wm-subgroups in terms of Wm-groups, and investigate some of their properties.

## 1. INTRODUCTION

Group theory and semi-group theory have developed in somewhat different directions in the past several decades. Group theory is the mathematical application of symmetry to an object to obtain knowledge of its physical properties. Group is the key part of it that acts in every area in which symmetry occurs. Lagrange, J. L. was usually credited with originating the theory of groups, which would become so important in 19th and 20th Century mathematics. Galois, E. work also laid the groundwork for further developments such as the beginnings of the field of abstract algebra, including area like group theory. Group theory is the tool that is used to determine symmetry and symmetry can help resolve many chemistry problems.

It is known that, a group G is an algebraic structure consisting of a non-empty set equipped with an operation on its elements that satisfies, associative law, identity law and inverse law and if the operation is commutative, then G is said to be a commutative group or an abelian group.

## 2. WM-GROUPS SETS

**Definition 2.1.** Let G be a non-empty set and * be a binary operation. Then, (G, ∗) is called a Wm-group if the following axioms hold:

(1) r * m * n = m * r * n, for every r, m, n ∈ G (semi commutativity).

(2) There is an element [r] in G such that r = r * [r] for each r ∈ G.

(3) There is an element $r^{-1}$ ∈ G such that [r] = $r^{-1}$ * r and $r^{-1}$ = $r^{-1}$ * [r] for each r ∈ G.

(4) There is a unique solution in G to the equation m * z = r, denoted by $m^{-1}$ * r, for every r, m ∈ G.

**Example 2.2.** Let G = {[r]= [m]= [$m^{-1}$], r = $r^{-1}$, m, $m^{-1}$} be a set and * be an operation on G defined as follows:

| * | [r]= [m]= [$m^{-1}$] | r = $r^{-1}$ | m | $m^{-1}$ |
|---|---|---|---|---|
| [r]= [m]= [$m^{-1}$] | [r]= [m]= [$m^{-1}$] | r | m | $m^{-1}$ |
| r = $r^{-1}$ | r = $r^{-1}$ | [r]= [m]= [$m^{-1}$] | $m^{-1}$ | m |
| m | m | $m^{-1}$ | r = $r^{-1}$ | [r]= [m]= [$m^{-1}$] |

Then, (G, ∗) is a Wm-group.

**Remark 2.3.** It is clear that every commutative group is Wm-group.
**Remark 2.4.** Wm-groups and groups are independent in general.

**Example 2.5.** Let A = {1, 2, 3} and let

$$\begin{pmatrix}1&2&3\\1&2&3\end{pmatrix}=1,\qquad \begin{pmatrix}1&2&3\\2&1&3\end{pmatrix}=(12)$$
$$\begin{pmatrix}1&2&3\\1&3&2\end{pmatrix}=(23),\qquad \begin{pmatrix}1&2&3\\3&2&1\end{pmatrix}=(13)$$
$$\begin{pmatrix}1&2&3\\2&3&1\end{pmatrix}=(123),\qquad \begin{pmatrix}1&2&3\\3&1&2\end{pmatrix}=(132)$$

Then, $S_3$ = {1, (12), (13), (23), (123), (132)} and thus ($S_3$, ∘) is a group, but not Wm-group because (12) ∘ (13) ∘ (23) = (13) ≠ (13) ∘ (12) ∘ (23) = (12).

**Example 2.6.** Let G be the set of all integer numbers and ∗ be an operation on G defined by r ∗ m = −(r + m) for r, m ∈ G. Then, (G, ∗) is a Wm-group, but not a group.

**Proposition 2.7.** Let (G, ∗) be a Wm-group. If n * r = n * m, then r = m, for every r, m, n ∈ G.
**Proof.** Since n * m = n * m implies that m = $n^{-1}$ * (m * n) by Definition 2.1 (4). But n * r = n * m, hence r = $n^{-1}$ * (m * n). Thus, r = m.

**Proposition 2.8.** Let (G, ∗) be a Wm-group. Then, r * ($m^{-1}$ * n) = ($m^{-1}$ * r) * n = $m^{-1}$ * (r * n) for every r, m, n ∈ G.
**Proof.** Let z = $m^{-1}$ * n and w = $m^{-1}$ * r. Then, m * z = n, m * w = r and m * w * n = r * n = r * m * z = m * r * z. By Proposition 2.7, w * n = r * z. Furthermore, m * r * z = m * r * ($m^{-1}$ * n) = r * m * ($m^{-1}$ * n) = r * n. Hence, r * z = $m^{-1}$ * (r * n).

**Proposition 2.9.** If (G, ∗) be a Wm-group, then,
(1) [r] * m = m for every r, m ∈ G.
(2) r * $r^{-1}$ = [r] for every r ∈ G.
(3) [$m^{-1}$] = [m] for every m ∈ G.

* Corresponding author

(4) $(r^{-1})^{-1} = r$ for every $r \in G$.
**Proof.** (1) By Proposition 2.8, we have $[r] * m = (r^{-1} * r) * m = r * (r^{-1} * m)$. Let $z = r^{-1} * m$, then $r * z = m$. Hence $[r] * m = r * z = m$.
(2) $r * r^{-1} = r * (r^{-1} * [r]) = [r]$.
(3) Since $[m^{-1}] = (m^{-1})^{-1} * m^{-1}$ implies that $m^{-1} * [m^{-1}] = m^{-1}$, but $m^{-1} * [m] = m^{-1}$, then $m^{-1} * [m^{-1}] = m^{-1} = m^{-1} * [m]$. By Proposition 2.7, $[m^{-1}] = [m]$.
(4) $r^{-1} * r = [r] = [r^{-1}]$. Thus, $r = (r^{-1})^{-1} * [r^{-1}] = (r^{-1})^{-1}$.

**Proposition 2.10.** Let $(G, *)$ be a Wm-group. Then, $(r * m)^{-1} * n = m^{-1} * (r^{-1} * n) = r^{-1} * (m^{-1} * n)$ for every $r, m, n \in G$.
**Proof.** Let $z = r^{-1} * n$. Then $r * z = n$ and $r * m * (m^{-1} * z) = r * z = n$. Hence $(r * m)^{-1} * n = m^{-1} * z = m^{-1} * (r^{-1} * n)$. Since $r * m * (m^{-1} * z) = m * r * (m^{-1} * z)$, then $r * (m^{-1} * z) = m^{-1} * n$. Thus, $r^{-1} * (m^{-1} * n) = m^{-1} * z = m^{-1} * (r^{-1} * n)$.

**Proposition 2.11.** Let $(G, *)$ be a Wm-group. Then for every $r, m \in G$, we have,
(1) $[r * m] = [m]$.
(2) $(r * m)^{-1} = r^{-1} * m^{-1}$.
(3) $[[m]] = [m]$.
(4) $[m^{-1}] = [m] = [m]^{-1}$
(5) $r * m = m * r$ if and only if $[r] = [m]$.
**Proof.** (1) Let $r, m \in G$. Since G is Wm-group, then there is $z_1$ in G such that $r^{-1} * z_1 = m$, that is, $z_1 = (r^{-1})^{-1} * m$ implies that $r * m = (r^{-1})^{-1} * m = z_1 \in G$. By Definition 2.1 (3) and Propositions 2.8 and 2.10, we have $[r * m] = (r * m)^{-1} * (r * m) = r * ((r * m)^{-1} * m) = r * (r^{-1} * (r^{-1} * r)) = r^{-1} * m = [m]$.
(2) By Definition 2.1 and Proposition 2.10, we have $(r * m)^{-1} = (r * m)^{-1} * [r * m] = (r * m)^{-1} * [m] = r^{-1} * (m^{-1} * [m]) = r^{-1} * m^{-1}$.
(3) $[[m]] = [m] * [m]^{-1}$, by Proposition 2.9 (1), $[[m]] = [m]^{-1}$ and $[[m]] = (m * m^{-1})^{-1} = m^{-1} * m = [m]$.
(4) $[m]^{-1} = [m]^{-1} * [[m]] = [m]^{-1} * [m] = [[m]] = [m]$. By Proposition 2.9 (3), $[m^{-1}] = [m] = [m]^{-1}$.
(5) If $r * m = m * r$, then $[r] = [m * r] = [r * m] = [m]$.

Conversely, if $[r] = [m]$, then by Definition 2.1 (2), we have $r * m = r * m * [m] = m * r * [m] = m * r * [r] = m * r$.

**Proposition 2.12.** Let $(G, *)$ be a Wm-group. Then, for every $r, z, w \in G$, we have:
(1) If $r * z = r * w$, then $z = w * [z]$.
(2) If $z = [z] * w$, then $r * z = r * w$.
(3) If $r^{-1} * z = r^{-1} * w$, then $z = w * [z]$.
(4) If $z = [z] * w$, then $r^{-1} * z = r^{-1} * w$.
**Proof.** We only prove (1) and (2), the other parts can be proved similarly.
(1) If $r * z = r * w$, then $r^{-1} * r * z = r^{-1} * r * w$ and $[r] * z = [r] * w$. Then, $[r] * z * [z] = [r] * w * [z]$ and by Proposition 2.9 (1), $z * [z] = w * [z]$. Thus, by Definition 2.1 (2), $z = w * [z]$.
(2) If $z = [z] * w$, then $r * z = r * [z] * w = r * w$.

**Definition 2.13.** Let G be Wm-group and $\varphi \neq S \subseteq G$, then S is called a Wm-subgroup of G if S is a Wm-group.

**Proposition 2.14.** Let G be Wm-group and $\varphi \neq S \subseteq G$. Then, S is a Wm-subgroup of G if and only if $m^{-1} * r \in S$, for every $r, m \in S$.
**Proof.** Let S be a Wm-subgroup of G. Then, there is $z_0$ in S such that $m * z_0 = r$ for every r, m in S. But $z_0$ is a solution to the same equation in G. Hence by Definition 2.1 (2), $z_0 = m^{-1} * r$.

Conversely, since * is semi commutativity on G, then * is also semi commutativity on S. Let $[m]^{-1}$, $r \in S$, for all $r, m \in S$. If we take $m = r$, then $[r] = r^{-1} * a \in S$ and if we take $r = m$, then $[m] = m^{-1} * m \in S$. Now for any $m \in S$ as $[m] \in S$, we have $m^{-1} = m^{-1} * [m] \in S$. Let $r, m \in S$, implies that $m^{-1} * r \in S$. Then, $m * (m^{-1} * r) = r$ and $m^{-1} * r$ is a solution in S to $m * z = r$. Since any other solution in S to $m * z = r$ is also a solution in G, so $m^{-1} * r$ is unique solution in S. Thus, S is a Wm-subgroup of G.

**Proposition 2.15.** The union of two Wm-subgroups of a Wm-group is Wm-subgroup if and only if one is contained in the other.
**Proof.** Suppose that $S_1$ and $S_2$ are two Wm-subgroups of a Wm-group G. If $S_1 \subseteq S_2$, then $S_1 \cup S_2 = S_2$ and if $S_2 \subseteq S_1$, then $S_1 \cup S_2 = S_1$. In either cases we get $(S_1 \cup S_2, *)$ is a Wm-subgroup of $(G, *)$.

Conversely, let $S_1 \nsubseteq S_2$ and $S_2 \nsubseteq S_1$, then there is an element $r \in S_1$, but $r \notin S_2$, and there is an element $m \in S_2$, but $m \notin S_1$. Now $r, m \in S_1 \cup S_2$, then by Proposition 2.14, $m^{-1} * r \in S_1 \cup S_2$, so either $m^{-1} * r \in S_1$ or $m^{-1} * r \in S_2$. If $m^{-1} * r \in S_1$, then $m^{-1} * r * r^{-1} = m^{-1} \in S_1$, but $(m^{-1})^{-1} \in S_1^{-1}$ implies that $m \in S_1$, which is not true. Again, if $m^{-1} * r \in S_2$, then $m * m^{-1} * r = r \in S_2$, which is not true. Therefore, $S_1 \subseteq S_2$ and $S_2 \subseteq S_1$.

**Definition 2.16.** Let G be a Wm-group, S a Wm-subgroup of G and $r \in G$. Then, $r * S = \{r * s: s \in S\}$ is called a left coset of S by r and $S * r = \{s * r: s \in S\}$ is called a right coset of S by r.

**Proposition 2.17.** Let G be a Wm-group and S a Wm-subgroup of G. Define a relation K on G as follows: rKm if and only if $r = s * m$ for some $s \in S$. Then, K is an equivalence relation on G whose equivalence classes are precisely the right cosets of S by the elements of G.
**Proof.** Let s be an element in S. Then, $s^{-1} * s = [s] \in S$. Let $r \in G$. Then, $r = [s] * r$ and rKr. Moreover, if rKm, then $r = s * m$ for some element s in S. Then $m = s^{-1} * r$. But $s^{-1} = s^{-1} * [s] \in S$. Thus mKr. Now, if rKm and mKn then $r = s * m$ and $m = g * n$ for $s, g \in S$. Then $r = s * g * n$ and we have rKn. The equivalence relation K on G partitions G into disjoint equivalence classes. Let Er be the equivalence class of all elements of G equivalent to r. If $z \in Er$, then $z = s * r$ for some $s \in S$ and hence $z \in S * r$.

Conversely, if z is any element of the right coset $S * r$, then $z = s * r$ for some $s \in S$ and $z \in Er$.

**Proposition 2.18.** Let G be a Wm-group and S a Wm-subgroup of G. The operation * defined on G|S as follows: $(S * r) * (S * m) = S * (r * m)$ so that G|S is a Wm-group under operation *.
**Proof.** First, we show that this operation is well-defined. Let $S * r = S * r'$ and $S * m = S * m'$. Then, $r = s * r'$ and $m = g * m'$ for $s, g \in S$. Then $r * m = s * r' * g * m' = s * g * r' * m'$ and we have $r * m \in S * r' * m'$, that is, $r * m$ and $r' * m'$ are in the same right coset. This shows that the induced operation is well defined.

Semi commutativity of * is immediate consequences of definition. If $S * r$ is any element of G|S, then $S * [r] = S * r^{-1} * r = S * r^{-1} * S * r = [S * r]$ and $S * [r] \in G|S$. Now for any $S * r \in G|S$ as $S * [r] \in G|S$, we have $S * r^{-1} * S * [r] = S * r^{-1} * [r] = S * r^{-1} \in G|S$. Finally, if $S * r$ and $S * m$ are elements of G|S, then $S * m * S * (m^{-1} * r) = S * (m * (m^{-1} * r)) = S * r$ and we have a solution in G|S to the equation $(S * m) * (S * z) = S * r$. This solution is unique, for if $(S * m) * (S * z) = S * r = (S * m) * (S * w)$, then $m * z = s * m * w$ for $s \in S$. Then $m * z = m * s * w$ and by Proposition 2.7, $z = s * w$, whence $S * z = S * w$.

**Definition 2.19.** Let $G_1$ and $G_2$ be two Wm-groups. A mapping $f: G_1 \rightarrow G_2$ is called a homomorphism if for every r, m in $G_1$, $f(r * m) = f(r) * f(m)$.
If f is injective, then f is called an isomorphism.

**Proposition 2.20.** Let f: $G_1 \rightarrow G_2$ be a homomorphism. Then, for every r, m ∈ G, we have:
(1) $f(m^{-1} * r) = f(m)^{-1} * f(r)$.
(2) $[f(r)] = f([r])$.
(3) $f(r^{-1}) = f(r)^{-1}$.
**Proof.** (1) $f(m) * f(m^{-1} * r) = f(m * (m^{-1} * r)) = f(l)$. Hence, $f(m^{-1} * r) = f(m)^{-1} * f(r)$.
(2) $f([r]) = f(r^{-1} * r) = f(r)^{-1} * f(r) = [f(r)]$.
(3) $f(r^{-1}) = f(r^{-1} * [r]) = f(r)^{-1} * f([r]) = f(r)^{-1} * [f(r)] = f(r)^{-1}$.

**Proposition 2.21**. Let f: $G_1 \rightarrow G_2$ be a homomorphism. Then,
(1) $f(S) = \{f(s) | s \in S\}$ is a Wm-subgroup of $G_2$ for any Wm-subgroup S of $G_1$.
(2) $f^{-1}(S') = \{z \in G | f(z) \in S'\}$ is a Wm-subgroup of $G_1$ for every Wm-subgroup S' of $G_2$.
**Proof.** (1) If f(s) and f(g) in f(S), then s, g ∈ S and $f(g)^{-1} * f(s) = f(g^{-1} * s) \in f(S)$.
(2) If z, w ∈ $f^{-1}(S')$, then f(z) ∈ S' and f(w) ∈ S', whence $f(w)^{-1} * f(z) = f(w^{-1} * z) \in S'$ and $w^{-1} * z \in f^{-1}(S')$.

**Definition 2.22.** Let G be a Wm-group and S a Wm-subgroup of G. Then, S is called a normal Wm-subgroup of G if [z] ∈ S, for every z ∈ G.

**Proposition 2.23.** Let G be a Wm-group, S a Wm-subgroup of G and r ∈ G. Let $r^{-1} * S * r = \{r^{-1} * s * r | s \in S\}$. Then, S is normal Wm-subgroup in G if and only if $r^{-1} * S * r \subseteq S$ for every r ∈ G.
**Proof.** Let S be normal in G, then $r^{-1} * s * r = s * r^{-1} * r = s * [r] \in S$ for any s ∈ S. Hence, $r^{-1} * S * r \subseteq S$.

Conversely, let r ∈ G, then $[r] = r^{-1} * r = r^{-1} * [s] * r \in r^{-1} * S * r$ for at least one s ∈ S. Therefore, [r] ∈ S.

**Corollary 2.24.** Let G be a Wm-group, S a Wm-subgroup of G and r ∈ G. Then, S is a normal Wm-subgroup in G if and only if $r * S * r^{-1} \subseteq S$ for every r ∈ G.
**Proof.** The proof is similar to Proposition 2.23.

**Definition 2.25.** Let $G_1$ and $G_2$ be two Wm-groups. If f: $G_1 \rightarrow G_2$ is a homomorphism, then the kernel of f is defined as $\{r \in G_1 | f(r) = f([r])\}$ and is denoted by K.

**Proposition 2.26**. Let f: $G_1 \rightarrow G_2$ be a homomorphism and K be a kernel of f. Then, K is a normal Wm-subgroup of $G_1$.
**Proof.** Let r, m ∈ K. Then, $f(m^{-1} * r) = f(m)^{-1} * f(r) = f([m])^{-1} * f([r]) = f([m]^{-1}) * f([r]) = f([m]^{-1} * [r]) = f([m] * [r]) = f([r]) = f([m^{-1} * r])$. Thus, K is a Wm-subgroup of G. Moreover, if z ∈ $G_1$, then by Proposition 2.11 (3), $f([z]) = f([[z]])$ and thus [z] ∈ K. Therefore, K is a normal Wm-subgroup of $G_1$.

**Proposition 2.27**. Let f: $G_1 \rightarrow G_2$ be a homomorphism and K be a kernel of f. Then,
(1) If S is a normal Wm-subgroup of G and K ⊆ S, then f(S) is a normal Wm-subgroup of $G_2$ and $f^{-1}(f(S)) = S$.
(2) If S' is a normal Wm-subgroup of $G_2$, then $f^{-1}(S')$ is a normal Wm-subgroup of $G_1$ such that $K \subseteq f^{-1}(S')$ and $f(f^{-1}(S')) = S'$.
**Proof.** (1) By Proposition 2.21 (1), f(S) is a Wm-subgroup of $G_2$. Let w ∈ $G_2$. Then w = f(z) for some z ∈ $G_1$ and $[w] = [f(z)] = f([z]) \in f(S)$, whence S is normal in $G_1$. Hence f(S) is a normal Wm-subgroup of $G_2$. Let $z \in f^{-1}(f(S))$, so f(z) ∈ f(S). Then f(z) = f(s) for some s ∈ S and $f(s)^{-1} * f(z) = f(z)^{-1} * f(z)$. Then, $f(s^{-1} * z) = f(s)^{-1} * f(z) = f(z)^{-1} * f(z) = [f(z)] = f([z]) = f([s^{-1} * z])$ implies that $s^{-1} * z \in K$. Since K ⊆ H, we have $s^{-1} * z = s_1 \in S$ and $z = s * s_1$, hence $f^{-1}(f(S)) \subseteq S$. On the other hand if z ∈ S, then f(z) ∈ f(S) and $z \in f^{-1}(f(S))$. Therefore, $f^{-1}(f(S)) = S$.
(2) By Proposition 2.21 (2), $f^{-1}(S')$ is a Wm-subgroup of $G_1$. Moreover, if z ∈ $G_1$, then $f([z]) = [f(z)] \in S'$. Thus $[z] \in f^{-1}(S')$ and $f^{-1}(S')$ is a normal in $G_1$. Let p ∈ K. Then, $f(p) = f([p]) = [f(p)] \in S'$ and $p \in f^{-1}(S')$. Thus, $K \subseteq f^{-1}(S')$. Let $w \in f(f^{-1}(S'))$. Then w = f(z), where $z \in f^{-1}(S')$. Thus, w ∈ S'. On the other hand for any w ∈ S' we have w = f(z) for some z ∈ $G_1$ and $z \in f^{-1}(S')$. Then, $w = f(z) \in f(f^{-1}(S'))$ and we have $f(f^{-1}(S')) = S'$.

**Proposition 2.28.** Let S and R be Wm-subgroups of a Wm-group G, then
(1) S ∩ R is a Wm-subgroup of G.
(2) If S is normal in G, then S ∩ R is normal in R.
(3) If S and R are normal in G, then S ∩ R is normal in G.
**Proof.** (1) Let z, w ∈ S ∩ R. Then, $w^{-1} * z \in S$ and $w^{-1} * z \in R$ and hence $w^{-1} * z \in S \cap R$.
(2) If S is normal in G, then for any p ∈ R, [p] ∈ R and [p] ∈ S. Thus, [p] ∈ S ∩ R.
(3) If both S and R are normal in G, then for any r ∈ G, [r] ∈ S, [r] ∈ R and thus [r] ∈ S ∩ R.

**Proposition 2.29.** If S and R are Wm-subgroups of a Wm-group G such that S ⊆ R and S is normal in G, then R is normal in G.
**Proof.** Let r ∈ G. Then, [r] ∈ S and [r] ∈ R.

**Definition 2.30.** Let S and R be Wm-subgroups of a Wm-group G. Then, S * R denotes the subset {s * p| s ∈ S, p ∈ R} of G.

**Proposition 2.31.** Let S and R be Wm-subgroups of a Wm-group G, then
(1) S * R is a Wm-subgroup of G.
(2) If S is normal in G, then S * R is normal in G.
(3) If S and R are normal in G, then S * R = R * S.
Proof. (1) Let z, w ∈ S * R. Then z = s * g and w = p * u, where s, p ∈ H and g, u ∈ R. Then, $w^{-1} * z = (p * u)^{-1} * s * g = p^{-1} * u^{-1} * s * g = (p^{-1} * s) * (u^{-1} * g) \in S * R$. Thus, S * R is a Wm-subgroup of G.
(2) Since R is non-empty, $p_1 \in R$ for some $p_1$ and thus $[p_1] \in R$, whence S ⊆ S * R. If S is normal in G, then by Proposition 2.29, S * R is normal in G.
(3) If S and R are normal in G, then by Corollary 2.24, for s * p ∈ S * R we have $s * p = s * p * [p] = p * s * [p] = p * s * p * p^{-1} = p * p * s * p^{-1} \in p * R * S * p^{-1} \subseteq R * S$, and similarly $p * s = p * s * [s] = s * p * [s] = s * p * s * s^{-1} = s * s * p * s^{-1} \in s * S * R * s^{-1} \subseteq S * R$, for p * s ∈ R * S. Therefore, S * R = R * S.

**Proposition 2.32.** [G] = {r ∈ G| r = [r]} is a normal Wm-subgroup of a Wm-group G, and [G] = {[r]| r ∈ G}.
**Proof.** Let z, w ∈ [G]. Then $w^{-1} * z = [w]^{-1} * [z] = [w] * [z] = [z] = [w^{-1} * z]$. Furthermore, if z ∈ G, then [z] = [[z]] and [z] ∈ [G], whence [G] is normal in G and [G] = {[r]| r ∈ G}.

**Remark 2.33.** The Wm-subgroup [G] = {r ∈ G| r = [r]} is said to be the trivial Wm-subgroup of a Wm-group G.

**Additional Reading**
We refer the reader to the books (Adhikari and Adhikari 2003, 2004; Clifford and Preston 1961) for further details.

## شێوازەکێ تازە یێ گروپێت ئاڵوگوری لاواز

**کورتیا لێکولینێ:**

ئارمانجا ڤێ ڤەکولینێ رابونە ب بێناسەکرن و خوێندنا شێوازەکێ تازە یێ گروپا کو دهێتە نیاسین ب Wm-group دگەل کریارەکا تەنیا پشتبەست کرن ببەڵگەنەوێستی ئاڵوگوری, بێلایەنێ راستێ, دژە پێچاوانەیێ چەپێ. هەروەسا دێ رابین ب پێشکێشکرنا بیرۆکەێن کوسیتا راستێ, کولکەی Wm-group, هومومورفیزم, ناڤەروك و نورمال Wm-subgroup یێت هاتینە پێناسەکرن ب Wm-group و لدوماهیێ دیفچونا چەند سالوخەتێت ڤێ گروبی دێهێتەکرن.

**خلاصة البحث:**

الهدف من هذه البحث هو تعرف و دراسة صنف جديد من الزمر معرفة بال Wm-group مع عملية ثنائية الوحيدة اعتماداً على البديهيات شبه تبديلية, محايد الايمن و المعكوس الايسر. وكذلك نحن نقدم مفاهيم كوسيت الايمن, القسمة Wm-group , التشاكل, النوات و Wm-subgroup الاعتيادي تحددها بال Wm-group واخيراً التحقيق في بعض خصائصها.